%% file: agt-4-10.tex
\documentclass{gtart_h}

\input agtout

\lognumber{10}
\volumenumber{4}
\volumeyear{2004}
\papernumber{10}
\published{24 March 2004}
\pagenumbers{151}{175}
\received{23 September 2002}
\revised{5 September 2003}
\accepted{27 January 2004}

\usepackage{amsmath,amssymb}
\usepackage{amsbsy}
\usepackage{euscript}
\newtheorem{Thm}{Th\'eor\`eme}
\newtheorem{Prop}{Proposition}
\newtheorem{Lem}{Lemme}
\newtheorem{Cor}{Corollaire}
\theoremstyle{definition}
\newtheorem*{Rq}{Remarque}
\newtheorem*{Rqs}{Remarques}

\newtheorem{Def}{D\'efinition}

\newtheorem{Not}{Notation}
\newtheorem{Cond}{Condition}

\begin{document}
\title{Sur la r\'ealisation des modules instables}
\covertitle{Sur la r\noexpand\'ealisation des modules instables}
\asciititle{Sur la realisation des modules instables}
\author{DongHua Jiang}
\address{LAGA, Institut Galil\'ee, Universit\'e Paris Nord\\93430 
Villetaneuse, France}
\asciiaddress{LAGA, Institut Galilee, Universite Paris Nord\\93430 
Villetaneuse, France}
\email{donghua.jiang@polytechnique.org}

\begin{abstract}
In this article, we give some conditions on the structure of an unstable 
module, which are satisfied whenever this module is the reduced cohomology 
of a space or a spectrum. First, we study the structure of the sub-modules 
of $\Sigma^s \tilde{H}^\ast(B(\Bbb{Z}/2)^{\oplus d}; \Bbb{Z}/2)$, i.e., the 
unstable modules whose nilpotent filtration has length 1.  Next, we 
generalise this result to unstable modules whose nilpotent filtration has 
a finite length, and which verify an additional condition. The result says 
that under certain hypotheses, the reduced cohomology of a space or a 
spectrum does not have arbitrary large gaps in its structure. This result 
is obtained by applying Adams' theorem on the Hopf invariant and the 
classification of the injective unstable modules.

This work was carried out under the direction of L. Schwartz.
%\end{abstract}

%\begin{abstract}
{\bf R\'esum\'e}\qua
Dans cet article, on donne des restrictions sur la structure d'un module 
instable, qui doivent \^etre v\'erifi\'ees pour que celui-ci soit la 
cohomologie r\'eduite d'un espace ou d'un spectre. On commence par une 
\'etude sur la structure des sous-modules de $\Sigma^s 
\tilde{H}^\ast(B(\Bbb{Z}/2)^{\oplus d}; \Bbb{Z}/2)$, i.e., les modules 
instables dont la filtration nilpotente est de longueur 1. Ensuite, on 
g\'en\'eralise le r\'esultat aux modules instables dont la filtration 
nilpotente est de longueur finie, et qui v\'erifient une condition 
suppl\'ementaire. Le r\'esultat dit que sous certaines hypoth\`eses, la 
cohomologie r\'eduite d'un espace ou d'un spectre ne contient pas de 
lacunes de longueur arbitrairement grande. Ce r\'esultat est obtenu par 
application du c\'el\`ebre th\'eor\`eme d'Adams sur l'invariant de Hopf et 
de la classification des modules instables injectifs.

Ce travail est effectu\'e sous la direction de L. Schwartz.
\end{abstract}

\asciiabstract{In this article, we give some conditions on the structure of an unstable 
module, which are satisfied whenever this module is the reduced
cohomology of a space or a spectrum. First, we study the structure of
the sub-modules of Sigma^sH^*(B(Z/2)^{oplus d};Z/2), i.e., the unstable
modules whose nilpotent filtration has length 1. Next, we generalise
this result to unstable modules whose nilpotent filtration has a
finite length, and which verify an additional condition. The result
says that under certain hypotheses, the reduced cohomology of a space
or a spectrum does not have arbitrary large gaps in its
structure. This result is obtained by applying Adams' theorem on the
Hopf invariant and the classification of the injective unstable
modules.

This work was carried out under the direction of L. Schwartz.

Resume 

Dans cet article, on donne des restrictions sur la structure d'un
module instable, qui doivent etre verifiees pour que celui-ci soit la
cohomologie reduite d'un espace ou d'un spectre. On commence par une
etude sur la structure des sous-modules de Sigma^sH^*(B(Z/2)^{oplus
d};Z/2), i.e., les modules instables dont la filtration nilpotente est
de longueur 1.  Ensuite, on generalise le resultat aux modules
instables dont la filtration nilpotente est de longueur finie, et qui
verifient une condition supplementaire. Le resultat dit que sous
certaines hypotheses, la cohomologie reduite d'un espace ou d'un
spectre ne contient pas de lacunes de longueur arbitrairement
grande. Ce resultat est obtenu par application du celebre theoreme
d'Adams sur l'invariant de Hopf et de la classification des modules
instables injectifs.

Ce travail est effectue sous la direction de L. Schwartz.}

\asciikeywords{Operations de Steenrod; module instable; theoreme d'Adams; la classification des modules instables injectifs}

\maketitle

\section{Introduction}

En topologie alg\'ebrique, pour distinguer les espaces, on introduit des invariants, tels que l'homologie, la cohomologie et les groupes d'homotopie des espaces. Nous nous int\'eressons dans cet article \`a la cohomologie r\'eduite des espaces en tant que module instable sur l'alg\`ebre de Steenrod. Nous consid\'erons d'abord le cas $p = 2$, les g\'en\'eralisations pour les nombres $p$ premiers impairs seront donn\'ees dans la derni\`ere section.

Un probl\`eme central sur les modules instables est de savoir quand un tel module est la cohomologie r\'eduite d'un espace. Un r\'esultat c\'el\`ebre de J.F. Adams impose des restrictions fortes \`a un module instable pour qu'il soit la cohomologie r\'eduite d'un espace. Voici le r\'esultat d'Adams dont il est question:

\begin{Thm}[Adams \cite{Ad60}] \label{Th1}
Soit $X$ un espace ou un spectre, $k \geq 4$, soit $x \in H^n(X; \Bbb{Z}/2)$ tel que $Sq^{2^i} x = 0$, $\forall$ $i<k$, alors $Sq^{2^k} x \in \sum_{i<k} {\rm Im}(Sq^{2^i})$. 
\end{Thm}

\begin{Def}
Un module sur l'alg\`ebre de Steenrod $M$ est un {\sl module instable} si pour tout \'el\'ement $x \in M$, $Sq^i x = 0$ quand $i > |x|$. Ici, $|x|$ d\'esigne le degr\'e de $x$.
\end{Def}

Comme $Sq^0$ est l'identit\'e, ceci implique que les modules instables sont triviaux en degr\'e strictement inf\'erieur \`a z\'ero.

\begin{Def}
Par {\sl lacune} de longueur $d$ dans un module instable $M$, on entend une suite d'entiers $I = \{ i, \cdots, i+d-1 \}$ telle que $M^j = \{ 0 \}$, si $j \in I$, $M^{i-1} \neq \{ 0 \}$, $M^{i+d} \neq \{ 0 \}$. On note cette lacune par $(i-1, i+d)$ ou $(i-1, i+d-1]$.
\end{Def}

Issue du th\'eor\`eme d'Adams, une question int\'eressante est de savoir si dans la cohomologie mod 2 d'un espace, il peut exister ou non des lacunes de longueur arbitrairement grande. Dans cet article, on d\'emontre que c'est impossible sous certaines hypoth\`eses suppl\'ementaires sur la structure du module instable.

Nous devons rappeler, pour \'enoncer ces conditions, diverses d\'efinitions. Rappelons qu'un module $M$ est connexe si $M^{\leq 0} = \{0\}$, un module instable est donc connexe si $M^0 = \{0\}$.

\begin{Def}
La {\sl suspension} d'un module instable $M$ est le module $\Sigma M$ tel que $(\Sigma M)^n = M^{n-1}$, $\forall$ $n$.
\end{Def}

\begin{Def}
Un module instable $M$ est {\sl r\'eduit} si le morphisme $Sq_0: M \to M$ d\'efini par $Sq_0(x) = Sq^{|x|}(x)$, $\forall$ $x \in M$, est injectif.
\end{Def}

On va se restreindre dans la suite \`a \'etudier des modules instables dont l'enveloppe injective est somme directe finie d'objets injectifs ind\'ecomposables. D'apr\`es la classification des ${\mathcal U}-$injectifs (Lannes-Schwartz, \cite{LS89}), on sait que pour un tel module instable r\'eduit $M$, il existe des entiers $d$ et $\alpha_d$ tels que $M$ se plonge dans $H^\ast (B(\Bbb{Z}/2)^{\oplus d}; \Bbb{Z}/2)^{\oplus \alpha_d}$. Si $M$ est connexe, on peut supposer $\alpha_d = 1$. Donc pour \'etablir une propri\'et\'e pour les modules instables r\'eduits, il suffit de le faire pour les sous-modules instables de $H^\ast (B(\Bbb{Z}/2)^{\oplus d}; \Bbb{Z}/2)^{\oplus \alpha_d}$. Dans la suite on supposera $\alpha_d = 1$, les d\'emonstrations s'\'etendent sans probl\`eme.

\begin{Def}
(Schwartz \cite{Sc94})\qua
Un module instable $M$ est {\sl $s-$nilpotent} s'il est l'union de ses sous-modules ayant une filtration finie dont les quotients sont des $s-$\`eme suspensions.
\end{Def}

Soit $\mathcal U$ la cat\'egorie des modules instables. On d\'esigne ${\mathcal N}il_s$ la sous-cat\'egorie ab\'elienne pleine de $\mathcal U$ des modules $s-$nilpotents. La sous-cat\'egorie ${\mathcal N}il_s$ est \'epaisse (voir \cite{Ga62}, \cite{Sc01}). On a une filtration de $\mathcal U$:
$$\cdots \subset {\mathcal N}il_2 \subset {\mathcal N}il_1 = {\mathcal N}il \subset {\mathcal N}il_0 = {\mathcal U}.$$

Soit $nil_s: {\mathcal U} \to {\mathcal N}il_s$ l'adjoint \`a droite de l'inclusion ${\mathcal N}il_s \hookrightarrow {\mathcal U}$, $nil_sM$ est le plus grand sous-module d'un module instable $M$ dans ${\mathcal N}il_s$ et on a la filtration nilpotente de $M$:
$$\cdots \subset nil_2M \subset nil_1M \subset nil_0M = M.$$

\begin{Prop}[\cite{Ku95}, \cite{Sc86}] Soit $M$ un module instable. Alors le quotient $nil_sM / nil_{s+1}M$ est la $s-$\`eme suspension d'un module instable r\'eduit $R_s$, donc
$$nil_sM / nil_{s+1}M \cong \Sigma^s R_s.$$ 
\end{Prop}

\begin{Def}
La filtration nilpotente d'un module instable $M$ est {\sl de longueur finie} s'il existe un $n \geq 0$ tel que $nil_n M =0$.
\end{Def}

\begin{Def}\label{De1}
Soit $M$ un module instable connexe r\'eduit non-trivial. On d\'esigne par $n_1 < n_2 < \cdots$ les degr\'es $n$ tels que $M^n \neq \{ 0 \}$. Supposons que $M$ se plonge dans $H^\ast (B(\Bbb{Z}/2)^{\oplus d}; \Bbb{Z}/2)$. Le module instable $M$ sera dit {\sl de type $\mathcal T$}, s'il contient une lacune $(s, s+l]$ avec $s \geq n_1$ et
$$l \geq \max\{ 2^{d+4}, n_{j+1}-n_j\ |\ j=1, \cdots, 1+(d-1)2^{d-2} \}.$$
\end{Def}

\begin{Rq}
Le module $M$ est n\'ecessairement infini car $M$ est r\'eduit non-trivial.
\end{Rq}

Le r\'esultat principal de cet article est le th\'eor\`eme suivant:

\begin{Thm}\label{Thm}
Soit $M$ un ${\mathcal A}_2-$module qui est une suspension it\'er\'ee d'un sous-module de type $\mathcal T$ de ${\tilde H}^\ast (B(\Bbb{Z}/2)^{\oplus d}; \Bbb{Z}/2)$. Alors $M$ n'est pas r\'ealisable, i.e., il n'existe aucun espace $X$ tel que $M = \tilde{H}^\ast(X; \Bbb{Z}/2)$.
\end{Thm}

En fait le th\'eor\`eme d'Adams s'applique aussi aux spectres. Il en est donc de m\^eme du th\'eor\`eme pr\'ec\'edent, la suspension it\'er\'ee peut \^etre positive ou n\'egative et le module n'est ni la cohomologie r\'eduite d'un espace ni la cohomologie d'un spectre.

Une g\'en\'eralisation de ce th\'eor\`eme est faite sous certaines hypoth\`eses pour les modules instables connexes ayant une filtration nilpotente de longueur finie.

\begin{Def}\label{De6}
Soit $M$ un module instable infini connexe dont la filtration nilpotente est de longueur finie. Les quotients $nil_sM / nil_{s+1}M$ non-triviaux s'\'ecrivent sous la forme $\Sigma^{m_i} R_{m_i}$, $R_{m_i}$ r\'eduits, $i=1, \cdots, t$, $m_1 < \cdots < m_t$. Notons que l'un au moins des $R_{m_i}$ est infini. Supposons qu'il existe des entiers $d$ et $\alpha_d$ tels que tous les $R_{m_i}$ se plongent dans $H^\ast (B(\Bbb{Z}/2)^{\oplus d}; \Bbb{Z}/2)^{\oplus \alpha_d}$. Notons $I \subset \{ 1, \cdots, t \}$ le sous-ensemble des $i$ tels que $R_{m_i}$ soit infini, et soit $n_{1,i} < n_{2,i} < \cdots$ les degr\'es en lesquels ce module est non-trivial.

Soit $\delta$ tel que $2^\delta \geq t > 2^{\delta-1}$. Le module instable $M$ sera dit {\sl de type $\mathcal T$} s'il contient une lacune $(s, s+l]$ avec $s \geq \min\{ m_i+n_{1,i}\ |\ i \in I \}$ et
$$l \geq \max\{ (m_t+1)2^{d+4}, n_{j+1,i}-n_{j,i}\ |\ i \in I, j=1, \cdots, 1+(d+\delta-1)2^{d-2} \}.$$
\end{Def}

\begin{Cond}\label{Cond1}
Soit $M$ un module instable connexe dont la filtration nilpotente est de longueur finie. En utilisant les notations introduites dans la d\'efinition pr\'ec\'edente, on dira que $M$ v\'erifie la condition \ref{Cond1} si
$$\begin{array}{cc}
m_{i+1}-m_i \neq 1,2,4,8,& 1 \leq i \leq t-1, \\
m_{i+2}-m_i \neq 8,& 1 \leq i \leq t-2,
\end{array}$$
c'est-\`a-dire, $m_j-m_i = 2^\beta$ n'a pas de solution pour $1 \leq i,j \leq t$ et $0 \leq \beta \leq 3$.
\end{Cond}

\begin{Thm}\label{Th4}
Soit $M$ un module qui est une suspension it\'er\'ee (positive ou n\'egative) d'un module instable connexe dont la filtration nilpotente est de longueur finie, qui est de type $\mathcal T$ et v\'erifie la condition \ref{Cond1}. Alors $M$ n'est pas r\'ealisable, i.e., il n'existe aucun espace ou spectre $X$ tel que $M = \tilde{H}^\ast(X; \Bbb{Z}/2)$.
\end{Thm}

\begin{Cor}
La longueur des lacunes ne peut pas \^etre arbitrairement grande dans un module instable connexe r\'ealisable dont la filtration nilpotente est de longueur finie et qui v\'erifie la condition \ref{Cond1}. \qed
\end{Cor}

Dans cet article, on ne consid\`ere que les modules dont l'enveloppe injective est somme directe finie de modules injectifs ind\'ecomposables. Les r\'esultats obtenus sont cons\'equences du th\'eor\`eme d'Adams et de la classification de Lannes-Schwartz.

Voici quelques d\'etails sur le plan de cet article. Dans la section 2, on d\'efinit des op\'erations $Q_t^s$, $s, t \geq 0$, qui g\'en\'eralisent les op\'erations de Milnor. La section 3 contient un r\'esultat combinatoire. En utilisant ce r\'esultat, le th\'eor\`eme \ref{Thm} est d\'emontr\'e dans la section 4. Ensuite, le th\'eor\`eme \ref{Th4} est d\'emontr\'e dans la section 5. La derni\`ere section contient des g\'en\'eralisations pour le cas $p$ premier impair. Il y a un appendice \`a la fin sur les op\'erations $Q_t^s$.

L'auteur tient \`a remercier le rapporteur pour ses remarques et ses conseils, qui l'ont aid\'e  \`a \'eviter bien des impr\'ecisions dans les d\'efinitions et d\'emonstrations.

\section{Les op\'erations $Q_t^s$, $s, t \geq 0$}

Dans cette section, on d\'efinit les op\'erations $Q_t^s$, $s, t \geq 0$ et on donne bri\`evement leurs propri\'et\'es utilis\'ees dans les sections suivantes. Pour plus de d\'etails sur ces op\'erations, on renvoie le lecteur \`a l'appendice.

\begin{Def} Les op\'erations $Q_t^s$, $s, t \geq 0$ sont d\'efinies r\'ecursivement comme suit:

(1)\qua $Q_0^s = Sq^{2^s}$;

(2)\qua $Q_{t+1}^s = [Sq^{2^{s+t+1}}, Q_t^s]$.
\end{Def}

\begin{Not}
On note souvent $Q_t^0$ par $Q_t$, qui est la notation usuelle de l'op\'eration de Milnor concern\'ee \cite{Mi58}.
\end{Not}

Pour \'etablir les propri\'et\'es de ces op\'erations $Q_t^s$, on a besoin d'introduire quelques notations.

\begin{Not}
Le symbole $(n_1, \cdots, n_d)$ d\'esignera le mon\^ome $u^{n_1} \otimes \cdots \otimes u^{n_d}$ ou $x_1^{n_1} \cdots x_d^{n_d}$ dans $H^\ast
(B(\Bbb{Z}/2)^{\oplus d}; \Bbb{Z}/2)$ qui s'identifie \`a ${\Bbb F}_2[u]^{\otimes d}$ ou ${\Bbb F}_2[x_1, \cdots, x_d]$, $u$ et les $x_i$
\'etant de degr\'e 1. Un tel mon\^ome sera dit {\sl basique}.
\end{Not}

\begin{Not}
Comme plus haut, $Sq_0$ d\'esigne l'op\'eration d\'efinie dans un module instable par $Sq_0 x = Sq^{|x|} x$. On a donc $Sq_0^s(n_1, \cdots, n_d)=(2^sn_1, \cdots, 2^sn_d)$, et ${\rm Im}(Sq_0^s)$ est l'ensemble des \'el\'ements $x = \sum_{i \in I} (2^sn(i)_1, \cdots, 2^sn(i)_d)$ de $H^\ast (B(\Bbb{Z}/2)^{\oplus d}; \Bbb{Z}/2)$. Ici, $I$ est un ensemble d'indices $i$ qui indexent des diff\'erents $d-$uplets $(n(i)_1, \cdots, n(i)_d)$, $n(i)_\alpha$ peut \^etre nul.
\end{Not}

\begin{Lem}\label{Le0}
Soit $M$ un module instable, on a pour tout $n \geq 1$,
$$Sq^{2n} Sq_0 x = Sq_0 Sq^{n} x, \qquad \forall\ x \in M.\eqno{\qed}$$
\end{Lem}

De la d\'efinition de $Q_t^s$ et de $Sq_0$, on d\'eduit que:

\begin{Lem}\label{Le1}
Soit $M$ un module instable, on a $\forall$ $x \in M$,
$$Q_t^{s+r}Sq_0^s x = Sq_0^sQ_t^r x,\qquad \forall\ r, s, t.$$
En particulier,
$$Q_t^sSq_0^s x = Sq_0^sQ_t x,\qquad \forall\ s, t.$$
\end{Lem}

\begin{Cor}\label{Pr1}
Soient $M, N$ deux modules instables, et soient $l, r, s, t \geq 0$.

{\rm(1)}\qua $\forall$ $x \in M$,
$Q_r^sQ_t^sSq_0^s x = Q_t^sQ_r^sSq_0^s x$ et $(Q_t^s)^2Sq_0^s x = 0$.

{\rm(2)}\qua $\forall$ $x \in Sq_0^s(M)$ et $y \in Sq_0^s(N)$,
$Q_t^s(x \otimes y) = Q_t^sx \otimes y + x \otimes Q_t^sy$.

{\rm(3)}\qua Soit $u$ le g\'en\'erateur de $H^\ast(B(\Bbb{Z}/2); \Bbb{Z}/2)$ en degr\'e 1,

$Q_t^sSq_0^su^{2l} = 0$ et $Q_t^sSq_0^su^{2l+1} = Sq_0^su^{2l+2^{t+1}} = u^{2^s(2l+2^{t+1})}$.
\end{Cor}

Le lemme \ref{Le0} est une cons\'equence directe de la d\'efinition de $Sq_0$, sa d\'emonstration est laiss\'ee au lecteur. Pour le lemme \ref{Le1} et le corollaire \ref{Pr1}, leurs d\'emonstrations se trouvent dans l'appendice.

\section{Un r\'esultat combinatoire}\label{Se3}

Dans cette section, on \'etablit d'abord un r\'esultat combinatoire. Ensuite, on l'applique \`a un \'el\'ement quelconque de $H^\ast (B(\Bbb{Z}/2)^{\oplus d}; \Bbb{Z}/2)$ pour obtenir des contraintes impos\'ees par certaines conditions d'annulation induites par l'existence de lacunes.

Soit un \'el\'ement $x \in H^\ast (B(\Bbb{Z}/2)^{\oplus d}; \Bbb{Z}/2)$, $x = \sum_{i \in I} (n(i)_1, \cdots, n(i)_d)$ est somme de mon\^omes basiques deux \`a deux distincts $(n(i)_1, \cdots, n(i)_d)$, $i \in I$.

Pour commencer, on d\'efinit quelques notations combinatoires.

\begin{Def}
Soit $g \geq 0$, on dira qu'il y a un {\sl $g-$\'echange} entre deux mon\^omes basiques $\alpha$ et $\beta$ s'il existe $i$ et $j$ tels que ces deux mon\^omes constituent, \`a un ordre (entre $i$ et $j$) pr\`es, une paire de la forme
$$\alpha = (u_1, \cdots, u_{i-1}, {2u_i+1}, u_{i+1}, \cdots, u_{j-1}, {2u_j+2^{g+1}}, u_{j+1}, \cdots, u_d)$$
$$\beta = (u_1, \cdots, u_{i-1}, {2u_i+2^{g+1}}, u_{i+1}, \cdots, u_{j-1}, {2u_j+1}, u_{j+1}, \cdots, u_d).\leqno{\hbox{et}}$$
On dira plus pr\'ecis\'ement, s'il y a lieu, qu'il y a un {\sl $g-$\'echange en $i-$\`eme position} pour le mon\^ome $\alpha$ avec le mon\^ome $\beta$.
\end{Def}

\begin{Rq}
C'est l'annulation sous l'action de l'op\'eration de Milnor $Q_g$ sur un \'el\'ement $x$ qui sugg\`ere cette d\'efinition, puisque $Q_g$ est une d\'erivation.
\end{Rq}

\begin{Def}
On dira qu'il y a une {\sl $(l, s)-$cha{\^\i}ne}, $l \leq s$, entre deux mon\^omes basiques $\alpha$ et $\beta$ d'un sous-ensemble de l'ensemble des mon\^omes basiques d'un \'el\'ement $x \in H^\ast (B(\Bbb{Z}/2)^{\oplus d}; \Bbb{Z}/2)$ s'il existe des mon\^omes basiques:
$$\alpha = \alpha_0,\ \alpha_1,\ \cdots,\ \alpha_t = \beta$$
dans ce sous-ensemble tels qu'il y ait un $m-$\'echange, $l \leq m \leq s$, entre $\alpha_i$ et $\alpha_{i+1}$ pour tout $i=0, \cdots, t-1$.
\end{Def}

\begin{Def}
Soit $x \in H^\ast (B(\Bbb{Z}/2)^{\oplus d}; \Bbb{Z}/2)$. On dira qu'un sous-ensemble $S$ de l'ensemble des mon\^omes basiques de $x$ admet $T \subset \{ 1, \cdots, d \}$ pour {\sl support}, si pour tout mon\^ome basique $x_1^{\alpha_1} \cdots x_d^{\alpha_d}$ appartenant \`a $S$ et pour tout $i \in T$, l'exposant $\alpha_i$ ne d\'epend que de $S$ et pas du m\^onome basique choisi et est de plus pair. On suppose de plus $T$ maximal parmi les sous-ensembles de $\{ 1, \cdots, d \}$ ayant cette propri\'et\'e.

On note $\tau = \#T$ que l'on appellera la {\sl taille} de $T$, les mon\^omes basiques de $S$ ont donc $\tau$ exposants en commun et s'\'ecrivent tous sous la forme $y^2 z$ o\`u $y$ d\'epend de $\tau$ variables $x_i$ et ne d\'epend pas du mon\^ome basique choisi; $z$ d\'epend lui de $d-\tau$ variables et du mon\^ome basique choisi.
\end{Def}

\begin{Def}
Un sous-ensemble de l'ensemble des mon\^omes basiques d'un \'el\'ement $x \in H^\ast (B(\Bbb{Z}/2)^{\oplus d}; \Bbb{Z}/2)$ est appel\'e une {\sl $(l, s)-$classe}, $l \leq s$, de support $T \subset \{ 1, \cdots, d \}$, si la condition suivante a lieu: pour tout mon\^ome basique $\alpha$ dans ce sous-ensemble, il existe au moins une position $i$ dont l'exposant est impair; pour toutes ces positions $i$ et tous les $m$, $l \leq m \leq s$, il existe un mon\^ome $\beta$ dans le sous-ensemble et un $m-$\'echange en $i-$\`eme position pour $\alpha$ avec $\beta$.
\end{Def}

\begin{Rq}
C'est l'annulation sous l'action des op\'erations $Q_m$, $l \leq m \leq s$, sur un \'el\'ement $x$ qui sugg\`ere cette d\'efinition, puisque les op\'erations $Q_m$ sont des d\'erivations.
\end{Rq}

Voici la propri\'et\'e fondamentale des $(l, s)-$classes de support $T$:

\begin{Prop}\label{Pr2}
Pour toute $(l, s)-$classe de support $T$ d'un \'el\'ement $x \in H^\ast (B(\Bbb{Z}/2)^{\oplus d};$ $\Bbb{Z}/2)$, on a $s-l+\#T \leq d-2$.
\end{Prop}

{\bf D\'emonstration}\qua Consid\'erons une $(l, s)-$classe admettant $T = T_s$ pour support, soit $\tau_s$ sa taille. Pour $1 \leq t \leq s-l$, on va construire r\'ecursivement des $(l, s-t)-$classes de support $T_{s-t}$ de taille $\tau_{s-t}$ telles que
$$\tau_{s-t} \geq \tau_{s-t+1}+1.$$

Pour $t = s-l$, on aura une $(l, l)-$classe de support $T_l$ dont la taille $\tau_l$ sera telle que $\tau_l \geq \tau_{l+1}+1 \geq \cdots \geq \tau_s+s-l = \#T+s-l$. Comme cette classe comporte des $l-$\'echanges, on a $\tau_l \leq d-2$. D'o\`u,
$$d-2 \geq \tau_l \geq \#T+s-l.$$

Supposons avoir construit une $(l, s-t+1)-$classe de support $T_{s-t+1}$ de taille $\tau_{s-t+1}$. On va construire une sous$-(l, s-t)-$classe, de la $(l, s-t+1)-$classe initiale, dont le support sera obtenu par adjonction \`a $T_{s-t+1}$ d'une position o\`u l'exposant d'un certain mon\^ome $\beta$ prend une valeur paire.

On consid\`ere parmi les exposants impairs qui apparaissent dans les mon\^omes basiques de la $(l, s-t+1)-$classe la valeur maximale, soit $2a+1$. Notons qu'il apparait n\'ecessairement des exposants impairs car il y a des $m-$\'echanges, $l \leq m \leq s-t+1$. On suppose que cet exposant appara{\^\i}t en position $p$ d'un mon\^ome basique $\alpha$ de la $(l, s-t+1)-$classe. Soit alors $\beta$ un mon\^ome dans la $(l, s-t+1)-$classe tel qu'il existe un $(s-t+1)-$\'echange en position $p$ pour $\alpha$ avec $\beta$. Le mon\^ome $\beta$ existe par hypoth\`ese.

Si on d\'esigne par $\psi_p$, l'exposant en position $p$ d'un mon\^ome basique $\psi$, on a alors $\alpha_p = 2a+1$ et $\beta_p = 2a+2^{s-t+2}$.

\begin{Lem}
Pour tout mon\^ome basique $\gamma$ d'une $(l, s-t)-$cha{\^\i}ne contenue dans la $(l, s-t+1)-$classe et contenant $\beta$ la valeur $\gamma_p$ de l'exposant en position $p$ est $2a + 2^{s-t+2}$.
\end{Lem}

\begin{proof}[D\'emonstration]
Raisonnons par l'absurde et choisissons une $(l, s-t)-$cha{\^\i}ne contenue dans la $(l, s-t+1)-$classe qui ne satisfasse pas \`a cette condition et soit de longueur minimale. Soit $\beta = \beta_0$, $\cdots$, $\beta_u = \gamma$ cette cha{\^\i}ne. L'exposant $\gamma_p$ est impair. Il y a un $m-$\'echange, $l \leq m \leq s-t$, entre $\beta_{u-1}$ et $\gamma$ en position $p$. Mais $(\beta_{u-1})_p = 2a+2^{s-t+2}$, donc $\gamma_p = 2a+2^{s-t+2}-2^{m+1}+1 > 2a+1$, en contradiction avec la maximalit\'e de $2a+1$.
\end{proof}

Consid\'erons alors l'ensemble $S$ des mon\^omes basiques de la $(l, s-t+1)-$classe tels qu'il existe une $(l, s-t)-$cha{\^\i}ne entre ces mon\^omes et $\beta$.

\begin{Lem}
L'ensemble $S$ est une $(l, s-t)-$classe de support contenant $T_{s-t+1} \cup \{p\}$.
\end{Lem}

\begin{proof}[D\'emonstration]
Comme chaque mon\^ome basique de l'ensemble $S$ est par d\'efinition un mon\^ome basique de la $(l, s-t+1)-$classe, il contient donc au moins un exposant impair. Pour montrer que $S$ est une $(l, s-t)-$classe, il faut encore montrer qu'en toute position $q$ o\`u un des mon\^omes de cet ensemble a un exposant impair, il y a pour tout $m$, $l \leq m \leq s-t$, un $m-$\'echange en position $q$ pour chacun de ces mon\^omes avec un autre mon\^ome dans $S$. Mais un tel mon\^ome existe par hypoth\`ese dans la $(l, s-t+1)-$classe et ce mon\^ome est alors par d\'efinition dans $S$ puisqu'il y a une $(l, s-t)-$cha{\^\i}ne \`a $\beta$. Clairement le suport $T_{s-t}$ contient $T_{s-t+1} \cup \{p\}$.
\end{proof}

Il reste \`a observer pourquoi on peut mener le processus jusqu'\`a $t = s-l$, car dans cette construction comme il y a des $m-$\'echanges, $l \leq m \leq s-t$, il y a des exposants impairs.

{\bf Fin de la d\'emonstration de la proposition} \qed

\begin{Cor}\label{Co1}
Soit $x$ un \'el\'ement de $H^\ast (B(\Bbb{Z}/2)^{\oplus d};$ $\Bbb{Z}/2)$ tel que $x \in {\rm Im}(Sq_0^s) - {\rm Im}(Sq_0^{s+1})$ et que $Q_t^sx = 0$, $\forall$ $p \leq t \leq q$. Alors on a $q-p \leq d-2$.
\end{Cor}

\begin{proof}[D\'emonstration]
Puisque $x = Sq_0^sx' \in {\rm Im}(Sq_0^s) - {\rm Im}(Sq_0^{s+1})$, il existe au moins un mon\^ome basique $\alpha$ de $x'$ avec au moins un exposant impair. L'ensemble des mon\^omes basiques de $x'$ qui sont dans une $(p, q)-$cha{\^\i}ne contenant $\alpha$ est une $(p, q)-$classe dont on note le support par $T$. Pr\'ecisons un peu. D'apr\`es la d\'efinition, chaque mon\^ome basique $\beta$ de cet ensemble contient au moins un exposant impair (\`a cause de l'existence d'un \'echange avec un autre mon\^ome basique de l'ensemble). On note par $I_\beta$ l'ensemble non vide des positions des exposants impairs dans $\beta$. Comme l'action de l'op\'eration $Q_t$ sur $\beta$ donne un mon\^ome qui contient un exposant pair en la m\^eme position, l'annulation de $Q_t^s$ sur $x$ entra{\^\i}ne l'existence d'un $m-$\'echange, $p \leq m \leq q$, en $i-$\`eme position, $i \in I_\beta$, pour $\beta$ avec un mon\^ome basique de l'ensemble.

Par cons\'equent, la proposition \ref{Pr2} nous donne $q-p \leq q-p+\#T \leq d-2$.
\end{proof}

\begin{Cor}\label{Co2}
Soit $x$ un \'el\'ement de $H^\ast (B(\Bbb{Z}/2)^{\oplus d};$ $\Bbb{Z}/2)$. Si ${\mathcal A}_2x$ contient une lacune $(|x|, |x|+l]$, $l \geq 2^k$ pour un certain $k \geq d-2$, alors $x \in {\rm Im}(Sq_0^{k-d+2})$.
\end{Cor}

\begin{proof}[D\'emonstration]
Soit $\alpha$ tel que $x = Sq_0^\alpha x' \in {\rm Im}(Sq_0^\alpha) - {\rm Im}(Sq_0^{\alpha+1})$. Si $\alpha > k$, on a $\alpha \geq k-d+2$. Si $\alpha \leq k$, alors l'existence de la lacune dans ${\mathcal A}_2x$ implique que $\forall$ $t=0, \cdots, k-\alpha$,
$$Sq_0^\alpha Sq^{2^t}x' = Sq^{2^{t+\alpha}}Sq_0^\alpha x' = Sq^{2^{t+\alpha}}x = 0.$$
Donc $Sq^{2^t}x' = 0$, $\forall$ $t=0, \cdots, k-\alpha$. Donc $Q_tx' = 0$, $\forall$ $t=0, \cdots, k-\alpha$. D'apr\`es le corollaire \ref{Co1}, on a $k-\alpha \leq d-2$, d'o\`u $\alpha \geq k-d+2$.
\end{proof}

\begin{Rqs}
(1)\qua Les \'enonc\'es de cette section sont aussi vrais pour un \'el\'ement quelconque de $H^\ast (B(\Bbb{Z}/2)^{\oplus d}; \Bbb{Z}/2)^{\oplus \alpha_d}$. Ci-dessus, on a trait\'e le cas o\`u $\alpha_d = 1$. Pour tenir compte du fait que l'on peut se placer dans $H^\ast (B(\Bbb{Z}/2)^{\oplus d}; \Bbb{Z}/2)^{\oplus \alpha_d}$, il faudrait compliquer un peu les notations en rajoutant un indice $1 \leq a \leq \alpha_d$. On dit que les $(n(i)_{1,a}, \cdots, n(i)_{d,a})$ sont les mon\^omes basiques de $x$.

(2)\qua Comme la suspension commute avec les op\'erations de Steenrod ($\sigma^q Sq^i = Sq^i \sigma^q$), on peut aussi \'etablir les \'enonc\'es similaires de ces corollaires pour une suspension quelconque de $H^\ast (B(\Bbb{Z}/2)^{\oplus d}; \Bbb{Z}/2)^{\oplus \alpha_d}$.
\end{Rqs}

\section{D\'emonstration du th\'eor\`eme \ref{Thm}}

Cette section est consacr\'ee \`a la d\'emonstration du th\'eor\`eme \ref{Thm}. Soit donc $M$ un module instable qui est la cohomologie r\'eduite d'un espace ou d'un spectre. Supposons de plus que $M$ est r\'eduit. Alors:

\begin{Thm}[Lannes-Schwartz \cite{LS89}] Un module instable r\'eduit (resp. r\'eduit et connexe) $M$ dont l'enveloppe injective est somme directe finie d'injectifs ind\'ecomposables est isomorphe \`a un sous-module de $H^\ast (B({\Bbb Z}/2)^{\oplus d}; {\Bbb Z}/2)^{\oplus \alpha_d}$ $(\alpha_d > 0)$ (resp. $H^\ast (B({\Bbb Z}/2)^{\oplus d}; {\Bbb Z}/2)$) pour $d$ assez grand.
\end{Thm}

{\bf D\'emonstration du Th\'eor\`eme \ref{Thm}}\qua
Dans la suite, on va d\'emontrer l'\'enonc\'e suivant: Soit $M$ un sous-module de type $\mathcal T$ de ${\tilde H}^\ast (B({\Bbb Z}/2)^{\oplus d}; {\Bbb Z}/2)$. Alors $M$ n'est pas r\'ealisable, i.e., il n'existe aucun espace ou spectre $X$ tel que $M = {\tilde H}^\ast (X; {\Bbb Z}/2)$. On note que, une fois cet \'enonc\'e est \'etabli, le th\'eor\`eme \ref{Thm} est aussi \'etabli.

On raisonne par l'absurde. Soit $M$ un sous-module de type $\mathcal T$ de ${\tilde H}^\ast (B(\Bbb{Z}/2)^{\oplus d}; \Bbb{Z}/2)$ qui est la cohomologie r\'eduite d'un espace ou d'un spectre. Reprenons les notations introduites avant le th\'eor\`eme \ref{Thm}: $M$ est non-trivial dans les degr\'es $n_1 < n_2 < \cdots$. Supposons que pour $n \geq n_1$, $(n, n+l]$ soit la premi\`ere lacune de longueur $l$ telle que:
$$l \geq \max\{ 2^{d+4}, n_{j+1}-n_j\ |\ j=1, \cdots, 1+(d-1)2^{d-2} \}.$$
Soit $k$ l'unique entier ($\geq d+4$) tel que $2^{k+1} > l \geq 2^k$.

Soit donc $x \in M^n$, $x \neq 0$. Alors ${\mathcal A}_2x$ contient une lacune $(n, n+l]$ et le corollaire \ref{Co2} entra{\^\i}ne $2^{k-d+2} | n$.

\begin{Lem}\label{Le3}
En degr\'e strictement inf\'erieur \`a $n$, il n'existe pas de degr\'es $m$ tels que $2^{k-d+2} {\not|} m$ et $M^m \neq \{0\}$.
\end{Lem}

{\bf D\'emonstration}\qua Supposons qu'en degr\'e strictement inf\'erieur \`a $n$, il existe des degr\'es $m$ tels que $2^{k-d+2} {\not|} m$ et $M^m \neq \{0\}$. Soit $m_0$ le plus grand de ces degr\'es, et soit $y \in M^{m_0}$, $y \neq 0$. On suppose que $y = Sq_0^\alpha z \in {\rm Im}(Sq_0^\alpha) - {\rm Im}(Sq_0^{\alpha+1})$. Comme $2^{k-d+2} {\not|} m_0$, on a $\alpha \leq k-d+1$.

Si ${\mathcal A}_2y$ contient une lacune $(|y|, n+l]$, le corollaire \ref{Co2} implique que $\alpha \geq k-d+2$, ce qui est impossible. Donc le plus bas degr\'e sup\'erieur ou \'egal \`a $m_0+1$, en lequel ${\mathcal A}_2y$ est non-trivial, est inf\'erieur \`a $n+l$ et est donc de la forme $2^{k-d+2}q =: p$ d'apr\`es l'hypoth\`ese de maximalit\'e de $m_0$.

Comme les $Sq^{2^h}$ engendrent multiplicativement ${\mathcal A}_2$, on a $p-m_0 =: 2^\beta$ (rappelons l'hypoth\`ese de minimalit\'e de $p$). En particulier, $Sq^{2^\beta}y$ est non nul. Comme $y \in {\rm Im}(Sq_0^\alpha)$, on a alors $\beta \geq \alpha$ car $Sq^{2^\beta} Sq_0^\alpha z = Sq_0^\alpha Sq^{2^{\beta-\alpha}} z$ qui est nul si $\alpha > \beta$. On va montrer que $\alpha \geq k-d$.

Supposons $\alpha \leq k-d$ et $t \geq 1$. Alors,

\begin{Lem}
$Q_t^\alpha y$ est nul, tant que son degr\'e est inf\'erieur ou \'egal \`a $n+l$.
\end{Lem}

\begin{proof}[D\'emonstration]
Par l'hypoth\`ese de maximalit\'e de $m_0$, on sait qu'il suffit de montrer que $2^{k-d+2}$ ne divise pas le degr\'e de $Q_t^\alpha y$. En effet on a:
$$\begin{array}{rcl}
|Q_t^\alpha y|&=& m_0+2^\alpha(2^{t+1}-1)\\
&=& 2^{k-d+2}q-2^\beta+2^{\alpha+t+1}-2^\alpha\\
&=& 2^{k-d+2}q+(2^{\alpha+t}+\cdots+2^{\alpha+1}+2^\alpha)-2^\beta.
\end{array}$$
Comme $t \geq 1$ et $\beta \geq \alpha$, on sait que ce degr\'e est un multiple impair de $2^\alpha$ pour $\beta > \alpha$ et que c'est un multiple impair de $2^{\alpha+1}$ pour $\beta = \alpha$. Donc $2^{\alpha+2}$ ne divise pas ce degr\'e. Puisque $\alpha+2 \leq k-d+2$, $2^{k-d+2}$ ne divise pas ce degr\'e non plus.
\end{proof}

Or pour $1 \leq t \leq k-\alpha-1$,
$$\begin{array}{rcl}
|Q_t^\alpha y|&=& m_0+2^\alpha(2^{t+1}-1)\\
&\leq& m_0+2^\alpha(2^{k-\alpha}-1)\\
&\leq& n+l.
\end{array}$$
Donc d'apr\`es le corollaire \ref{Co1}, $k-\alpha-2 \leq d-2$, alors $\alpha \geq k-d$ et par cons\'equent, $\beta \geq \alpha \geq k-d$.

D'apr\`es la d\'efinition de $\beta$, $p = m_0 + 2^\beta$, et par hypoth\`ese, $2^{k-d+2} {\not|} m_0$ et $2^{k-d+2} | p$, on a donc $\beta \leq k-d+1$, ce qui implique que $M$ contient une lacune $(m_0, p)$. Car sinon, il existe $p' \in (m_0, p)$ tel que $M^{p'} \neq \{ 0 \}$. Par l'hypoth\`ese de maximalit\'e de $m_0$, $2^{k-d+2} | p'$. Donc $2^\beta = p-m_0 > p-p' \geq 2^{k-d+2}$, ce qui est contradictoire au fait que $\beta \leq k-d+1$. L'existence de la lacune $(m_0, p)$ dans $M$ et le th\'eor\`eme d'Adams impliquent que $\beta \leq 3$. Or comme $\beta \geq k-d$ et $k \geq d+4$, on a $\beta \geq k-d \geq 4$, ceci implique qu'un tel $m_0$ n'existe pas.

{\bf Fin de la d\'emonstration du Lemme \ref{Le3}} \qed

Notons donc les degr\'es plus petits que $n$ pour lesquels $M$ est non-trivial comme suit
$$n = r_0 > r_1 > \cdots \quad \text{et}\quad \forall\ i,\ 2^{k-d+2} | r_i.$$
On a

\begin{Lem}
Pour tout $x_i$ de degr\'e $r_i$, ${\mathcal A}_2x_i$ contient la lacune $(r_i, n+l]$.
\end{Lem}

\begin{proof}[D\'emonstration]
On raisonne par l'absurde. Si l'\'enonc\'e est faux, on choisit un \'el\'ement $x_i$ de degr\'e maximal tel que ${\mathcal A}_2x_i$ contienne des \'el\'ements non nuls de degr\'e sup\'erieur \`a $r_i$ et inf\'erieur \`a $n+l$. On choisit dans $(r_i, n+l]$ le plus bas degr\'e en lequel ${\mathcal A}_2x_i$ est non-trivial. En utilisant la base multiplicative de ${\mathcal A}_2$, il est de la forme $r_i+2^\gamma$, $\gamma \geq 0$. Comme les degr\'es dans l'intervalle $(r_i, n+l]$ o\`u il y a des \'el\'ements non nuls sont divisibles par $2^{k-d+2}$, on a $\gamma \geq k-d+2 \geq 6$. Par cons\'equent, on a une lacune $(r_i, r_i+2^\gamma)$, $\gamma \geq 6$, dans ${\mathcal A}_2x_i$ avec $Sq^{2^\gamma}x_i \neq 0$ en degr\'e $r_i+2^\gamma$. Par la maximalit\'e de $x_i$, il n'y a aucun \'el\'ement $y$ de degr\'e sup\'erieur \`a $r_i$ tel que ${\mathcal A}_2y$ contienne des \'el\'ements non nuls en degr\'e sup\'erieur \`a $|y|$ et inf\'erieur \`a $n+l$. Donc l'\'el\'ement non nul $Sq^{2^\gamma}x_i$ ne peut pas \^etre dans $\sum_{j<\gamma} {\rm Im}(Sq^{2^j})$. Alors l'existence de cette lacune $(r_i, r_i+2^\gamma)$ dans ${\mathcal A}_2x_i$ est impossible \`a cause du th\'eor\`eme d'Adams.
\end{proof}

On montre alors par r\'ecurrence que:
$$\forall\ i \geq j2^{d-2},\qquad 2^{k-d+j+2} | r_i.$$
Le cas $j=0$ est d\'emontr\'e ci-dessus. Supposons que c'est vrai pour $j$, alors
$$\begin{array}{rcl}
n+l-r_{(j+1)2^{d-2}}&=& l+\sum_{i=0}^{(j+1)2^{d-2}-1}(r_i-r_{i+1})\\
&=&     l+\sum_{h=0}^{j}\sum_{i=h2^{d-2}}^{(h+1)2^{d-2}-1}(r_i-r_{i+1})\\
&\geq&  2^k+\sum_{h=0}^{j}2^{d-2}2^{k-d+h+2}\\
&=&     2^k+\sum_{h=0}^{j}2^{k+h}\\
&=&     2^{k+j+1}.
\end{array}$$
Puisque ${\mathcal A}_2 x_i$ contient une lacune $(r_i, n+l]$ et d'apr\`es le corollaire \ref{Co2}, on a donc $2^{k-d+j+3} | r_i$, $\forall$ $i \geq (j+1)2^{d-2}$.

Rappelons que par l'hypoth\`ese sur $l$, il y a bien (au moins) $(d-1)2^{d-2}+1$ valeurs pour l'indice $i$ de $r_i$. On peut donc poser $w = (d-1)2^{d-2}$, alors $2^{k+1} | r_w$ et $2^{k+1} | r_{w+1}$. L'intervalle $(r_{w+1}, r_w)$ est donc aussi une lacune de $M$, avec $r_{w+1} \geq n_1$ et $r_w \leq n$, d'une longueur $h$ telle que
$$\begin{array}{rcl}
h& \geq& 2^{k+1}-1 \\
& \geq& l \\
& \geq& \max\{ 2^{d+4}, n_{j+1}-n_j\ |\ j=1, \cdots, 1+(d-1)2^{d-2} \}.
\end{array}$$
Ceci est contradictoire au choix de $(n, n+l]$.

{\bf Fin de la d\'emonstration du Th\'eor\`eme \ref{Thm}} \qed

\section{D\'emonstration du th\'eor\`eme \ref{Th4}}

Dans cette section, on va \'etudier des modules instables dont la filtration nilpotente est de longueur finie. Un exemple trivial d'un tel module instable est un module instable quelconque de dimension finie. Un autre exemple est la suspension d'un module instable r\'eduit. La cohomologie d'un groupe fini ou du classifiant d'un groupe compact v\'erifie aussi cette hypoth\`ese \cite{HLS}.

Ci-dessous, on d\'emontre un r\'esultat sur la non-existence de grandes lacunes dans les modules instables connexes r\'ealisables dont la filtration nilpotente est de longueur finie qui v\'erifie la condition \ref{Cond1}.

{\bf D\'emonstration du Th\'eor\`eme \ref{Th4}}\qua
Comme dans la d\'emonstration du th\'eor\`eme \ref{Thm}, il suffit de prouver l'\'enonc\'e pour les modules instables connexes dont la filtration nilpotente est de longueur finie, qui est de type $\mathcal T$ et v\'erifie la condition \ref{Cond1}.

On raisonne par l'absurde. Soit donc $M$ un module instable connexe qui est la cohomologie r\'eduite d'un espace ou d'un spectre. Reprenons les notations introduites avant le th\'eor\`eme \ref{Th4}: les quotients $nil_sM / nil_{s+1}M$ non-triviaux s'\'ecrivent sous la forme $\Sigma^{m_i} R_{m_i}$, $R_{m_i}$ r\'eduits, $i=1, \cdots, t$, $m_1 < \cdots < m_t$. Tous les $R_{m_i}$, $i \in I$, soient non-triviaux dans les degr\'es $n_{1,i} < n_{2,i} < \cdots$. Soit $\delta$ tel que $2^\delta \geq t > 2^{\delta-1}$. Supposons que pour $n \geq \min\{ m_i+n_{1,i}\ |\ i \in I \}$, $(n, n+l]$ soit la premi\`ere lacune dans $M$ de longueur
$$l \geq \max\{ (m_t+1)2^{d+4}, n_{j+1,i}-n_{j,i}\ |\ i \in I, j=1, \cdots, 1+(d+\delta-1)2^{d-2} \}.$$
Soit $k$ l'unique entier ($\geq d+4$) tel que $2^{k+1} > l \geq 2^k.$

\begin{Lem}\label{Prop03}
Pour tout $x$ tel que $|x| \leq n$, le module ${\mathcal A}_2 x$ contient la lacune $(|x|, n+l]$.
\end{Lem}

{\bf D\'emonstration}\qua Pour montrer cela, on raisonne par l'absurde. A tout \'el\'ement $x \in M$, on associe son degr\'e de nilpotence, c'est-\`a-dire, l'entier $m_x$ tel que $x \in nil_{m_x}M - nil_{m_x+1}M$.

Soit $x$ non nul de degr\'e maximal tel que ${\mathcal A}_2x$ n'est pas r\'eduit \`a $\{0\}$ dans l'intervalle $(|x|, n]$. (${\mathcal A}_2x$ contient la lacune $(n, n+l]$ par hypoth\`ese.)

Soit donc $y \in {\mathcal A}_2x$ de degr\'e minimal tel que $y \neq 0$, $|x| < |y| \leq n$, $\bar{y} \in \Sigma^{m_y} R_{m_y}$ sa r\'eduction que l'on note $\sigma^{m_y}v$, $v \in R_{m_y}$. Alors comme ${\mathcal A}_2y$ contient une lacune $(|y|, n+l]$ (rappelons l'hypoth\`ese de maximalit\'e de $x$), le corollaire \ref{Co2} implique que $v \in {\rm Im}(Sq_0^{k-d+2})$. On sait donc que $|y|-m_y$ est divisible par $2^{k-d+2}$ et on \'ecrit $|y|-m_y =: 2^{k-d+2}l_y$.

Soit de m\^eme la r\'eduction $\bar{x} \in \Sigma^{m_x}R_{m_x}$. Notons $\bar{x}=\sigma^{m_x}u$, $u \in {\rm Im}(Sq_0^s) - {\rm Im}(Sq_0^{s+1})$.

\begin{Lem}
Pour $s \leq k-d$, $t \geq 1$, $2^{k-d+2}$ ne divise pas le degr\'e de $Q_t^s u$.
\end{Lem}

\begin{proof}[D\'emonstration]
En effet on a $|Q_t^s u| = |u|+2^s(2^{t+1}-1)$. Consid\'erons le ${\mathcal A}_2-$module engendr\'e par $u$ dans $R_{m_x}$. Si ${\mathcal A}_2 u$ contient la lacune $(|u|, n+l-m_x]$, le degr\'e de $u$ est divisible par $2^{k-d+2}$, et dans ce cas $2^{k-d+2}$ ne divise pas $|Q_t^s u|$ pour $s \leq k-d$.

Supposons que ${\mathcal A}_2 u$ ne contienne pas la lacune $(|u|, n+l-m_x]$ et soit $t = Sq^{2^\beta} u$ l'\'el\'ement non nul du plus bas degr\'e avec $|t| \leq n-m_x$. Comme ${\mathcal A}_2 t$ contient la lacune $(|t|, n+l-m_x]$, $2^{k-d+2} |\ |t|$. Pour la m\^eme raison que dans la d\'emonstration du lemme \ref{Le3}, on a $\beta \geq s$. Alors
$$\begin{array}{rcl}
|Q_t^s u|& =& |t|-2^\beta+2^s(2^{t+1}-1) \\
& =& 2^{k-d+2}q-2^\beta+2^{s+t+1}-2^s
\end{array}$$
et $2^{k-d+2}$ ne divise pas $|Q_t^s u|$.
\end{proof}

Supposons d'abord que $s \leq k-d$. Tant que $|Q_t^s x| \leq n+l$, on a donc n\'ecessairement $Q_t^s u = 0$ pour des raisons de degr\'e. En effet si $Q_t^s u \neq 0$, $\sigma^{m_x} Q_t^s u = Q_t^s \bar{x} = \overline{Q_t^s x}$ la r\'eduction de $Q_t^s x$ dans $\Sigma^{m_x} R_{m_x}$ dont le degr\'e de nilpotence est $m_x$ (qui est {\it a priori} plus grand que ou \'egal \`a celui de $x$, voir \cite{Sc86},\cite{Sc01}), et dont le degr\'e est de la forme (en appliquant le corollaire \ref{Co2} \`a la lacune $(|Q_t^s u|, n+l-m_x]$ dans ${\mathcal A}_2 (Q_t^s u)$)
$$m_x + 2^{k-d+2}f, \quad f \geq 0.$$
Donc pour que $Q_t^s u$ soit non nul, il faudrait que son degr\'e soit multiple de $2^{k-d+2}$.

Pour tout $t$ tel que $1 \leq t \leq k-s-1$,
$$\begin{array}{rcl}
|Q_t^s x|&=& |x|+2^s(2^{t+1}-1)\\
&\leq& |x|+2^s(2^{k-s}-1)\\
&\leq& n+l.
\end{array}$$
Donc d'apr\`es le corollaire \ref{Co1},
$$k-s-2 \leq d-2,\quad \text{soit}\ s \geq k-d$$
et donc $|x|-m_x =: 2^{k-d}l_x$.

Revenons alors \`a l'\'el\'ement $y$ non nul du plus bas degr\'e, sup\'erieur ou \'egal \`a $|x|+1$ dans ${\mathcal A}_2x$. Il est de la forme $Sq^{2^\alpha}x$, le th\'eor\`eme d'Adams implique que $\alpha \leq 3$. En effet l'hypoth\`ese de maximalit\'e de $x$ implique que pour tout \'el\'ement $z$ non nul dont le degr\'e est entre $|x|+1$ et $n$, ${\mathcal A}_2z$ est r\'eduit \`a $\{0\}$ dans l'intervalle $(|z|, n]$. Donc $Sq^{2^\alpha}x \not\in \sum_{i<\alpha} {\rm Im}(Sq^{2^i})$.

On a alors $|x| + 2^\alpha = |y|$, $y \in nil_{m_y}M$ et $m_y \geq m_x$. D'o\`u, $m_x + 2^{k-d}l_x + 2^\alpha = m_y + 2^{k-d+2}l_y$. Donc $m_x + 2^\alpha = m_y$ mod $2^{k-d}$. D'autre part, $2^{k+1} > l \geq (m_t+1)2^{d+4}$, ce qui implique que
$$\begin{array}{rcl}
2^{k-d}& >& 8(m_t+1)\\
& =& (m_x+m_y+2^\alpha) + (m_t-m_x) + (m_t-m_y) + (8-2^\alpha) + 6m_t\\
& \geq& m_x+m_y+2^\alpha.
\end{array}$$

$$|m_x+2^\alpha-m_y| \leq m_x+m_y+2^\alpha < 2^{k-d},$$
donc $m_x+2^\alpha = m_y$. Or cette \'egalit\'e n'a pas de solution \`a cause de la condition \ref{Cond1}, l'existence d'un tel $x$ est donc contradictoire. Donc pour tout \'el\'ement $x$ de $M$ de degr\'e inf\'erieur \`a $n$, ${\mathcal A}_2x$ contient une lacune $(|x|, n+l]$.

{\bf Fin de la d\'emonstration du Lemme \ref{Prop03}} \qed

On ach\`eve la d\'emonstration en appliquant la d\'emonstration du th\'eor\`eme \ref{Thm} \`a chaque $R_{m_i}$, $i \in I$. Plus pr\'ecis\'ement, on applique la (derni\`ere) partie de la d\'emonstration du th\'eor\`eme \ref{Thm} - concernant une r\'ecurrence sur la divisibilit\'e par une puissance de 2 des degr\'es inf\'erieurs \`a $n+l$ pour lesquels le module est non-trivial - \`a $R_{m_i}$ ($i \in I$). On veut donc obtenir des informations sur la divisibilit\'e par une puissance de 2 de ses degr\'es inf\'erieurs \`a $n+l$ auxquels on a soustrait $m_i$, pour lesquels $R_{m_i}$ est non-trivial, afin de montrer qu'on aboutit \`a une contradiction.

En effet, si on note $\forall$ $i \in I$,
$$(n-m_i \geq)\ r_{0,i} > r_{1,i} > \cdots$$
les degr\'es inf\'erieurs \`a $n+l-m_i$, pour lesquels $R_{m_i}$ est non-trivial, on a
$$\forall\ w \geq (d+\delta-1)2^{d-2},\qquad 2^{k+\delta+1} | r_{w,i}.$$
Donc en degr\'e inf\'erieur ou \'egal \`a $r_{w,i}$, $R_{m_i}$ ($i \in I$) ne contient que des lacunes de longueur plus grande que ou \'egale \`a $2^{k+\delta+1}-1$. Donc il existe un $d_0 \leq n$ tel qu'en degr\'e inf\'erieur ou \'egal \`a $d_0$, les $\Sigma^{m_i} R_{m_i}$ ($i \in I$) ne contiennent que des lacunes de longueur plus grande que ou \'egale \`a $2^{k+\delta+1}-1$, et on suppose de plus qu'il en existe au moins une en degr\'e inf\'erieur ou \'egal \`a $d_0$.

Maintenant on choisit une lacune de la plus petite longueur parmi toutes celles en degr\'e inf\'erieur ou \'egal \`a $d_0$ contenues dans l'un des $\Sigma^{m_i} R_{m_i}$, $i \in I$. Par le choix de cette lacune, disons $\Sigma^{m_{i_0}} R_{m_{i_0}}$ ($i_0 \in I$), on sait qu'elle ne contient aucun degr\'e en lequel $\Sigma^{m_{i_0}} R_{m_{i_0}}$ est non-trivial, et qu'elle contient au plus un degr\'e en lequel $\Sigma^{m_i} R_{m_i}$, $1 \leq i \leq t$, $i \neq i_0$, est non-trivial. Car sinon, elle contiendrait une lacune d'un des $\Sigma^{m_i} R_{m_i}$, $i \neq i_0$, en contradiction avec l'hypoth\`ese de minimalit\'e sur la longueur de la lacune choisie. Comme cette lacune contient au plus $t-1$ degr\'es en lesquels $M$ est non-trivial, il existe donc une lacune $(n', n'+l']$, avec $n' \geq \min\{ m_i+n_{1,i}\ |\ i \in I \}$ et $n'+l' \leq d_0$ ($\leq n$), de longueur
$$\begin{array}{rcl}
l'& \geq& \frac{1}{t} \cdot 2^{k+\delta+1} - 1 \\
& \geq& 2^{k+1} - 1 \\
& \geq& l \\
& \geq& \max\{ (m_t+1)2^{d+4}, n_{j+1,i}-n_{j,i}\ |\ i \in I, \\
& & \qquad \qquad \qquad j=1, \cdots, 1+(d+\delta-1)2^{d-2} \}.
\end{array}$$
Ceci est contradictoire au choix de $(n, n+l]$.

{\bf Fin de la d\'emonstration du Th\'eor\`eme \ref{Th4}}
\qed

\section{Le cas $p$ premier impair}

On indique bri\`evement les r\'esultats pour le cas $p$ premier impair. On donne d'abord les ingr\'edients essentiels, i.e., les op\'erations $Q_t^s$, $s, t \geq 0$, et $P_0^s$, $s \geq 0$. En tenant compte des signes, les formules et les r\'esultats combinatoires sont \'etablis de la m\^eme mani\`ere.

\begin{Def}
Les op\'erations $Q_t^s$, $s, t \geq 0$ sont d\'efinies r\'ecursivement comme suit:

(1)\qua $Q_0^0 = \beta$ et $Q_{t+1}^0 = [P^{p^t}, Q_t^0]$, ce sont des $Q_t$ d\'efinis par Milnor \cite{Mi58};

(2)\qua $Q_0^s = P^{p^{s-1}}$ et $Q_{t+1}^s = [P^{p^{s+t}}, Q_t^s]$, $s \geq 1$.
\end{Def}

\begin{Not}
Le symbole $(n_1, \cdots, n_d) = (2m_1+\epsilon_1, \cdots, 2m_d+\epsilon_d)$, $\epsilon_i =$ 0 ou 1 ($1 \leq i \leq d$), d\'esigne l'\'el\'ement $t^{\epsilon_1}u^{m_1} \otimes \cdots \otimes t^{\epsilon_d}u^{m_d} \in H^\ast (B(\Bbb{Z}/p)^{\oplus d}; \Bbb{Z}/p)$, $t$ \'etant de degr\'e 1 et $u$ \'etant de degr\'e 2. Un tel \'el\'ement est dit {\sl basique}.
\end{Not}

\begin{Not}
Comme d'habitude, $P_0$ d\'esigne l'op\'eration d\'efinie dans un module instable par
$$P_0 x = \left\{ \begin{array}{ll}
P^{|x|/2} x, & \text{ si } |x|=0 \text{ mod } 2;\\
\beta P^{(|x|-1)/2} x, & \text{ si } |x|=1 \text{ mod } 2.
\end{array} \right.$$
On a donc
$$P_0^s(2m_1+\epsilon_1, \cdots, 2m_d+\epsilon_d) = (2p^sm_1+2p^{s-1}\epsilon_1, \cdots, 2p^sm_d+2p^{s-1}\epsilon_d),$$
et ${\rm Im}(P_0^s)$ ($s \geq 1$) est l'ensemble des \'el\'ements $x = \sum_{j \in J} (2p^{s-1}l(j)_1, \cdots,$ $2p^{s-1}l(j)_d)$ de $H^\ast (B(\Bbb{Z}/p)^{\oplus d}; \Bbb{Z}/p)$. Par convention, $P_0^0 = id$. Ici, $J$ est un ensemble d'indices $j$ qui indexent des diff\'erents $d-$uplets $(l(j)_1, \cdots, l(j)_d)$, $l(j)_1, \cdots, l(j)_d = 0, 1$ mod $p$, $l(j)_\alpha$ peut \^etre nul.
\end{Not}

\begin{Rq}
Dans toute cette section, le symbole $P_0^s$ d\'esigne l'op\'eration $(P_0)^s$, qui est \'evidemment distincte de l'op\'eration de Milnor (utilis\'ee dans l'appendice).
\end{Rq}

Puisqu'on a pour tout \'el\'ement $x$ d'un module instable et pour tout $n \geq 1$,
$$P^{pn}P_0 x = P_0P^n x\qquad \text{ et }\qquad P^1P_0 x = P_0\beta x,$$
on peut donc \'etablir les propri\'et\'es de $Q_t^s$, $s, t \geq 0$, sur ${\rm Im}(P_0^s)$ \`a partir de celles de $Q_t = Q_t^0$. Une autre fa{\c{c}}on d'\'etablir les propri\'et\'es de $Q_t^s$ est d'utiliser le fait que
$$\psi^\ast(P^{p^s}) = P^{p^s} \otimes 1 + P^{p^s-1} \otimes P^1 + \cdots + 1 \otimes P^{p^s}$$
devient $P^{p^s} \otimes 1 + 1 \otimes P^{p^s}$, une d\'erivation sur ${\rm Im}(P_0^{s+1})$.

\begin{Lem}\label{Le5}
Soit $x$ un \'el\'ement de $H^\ast (B(\Bbb{Z}/p)^{\oplus d}; \Bbb{Z}/p)$ tel que $x \in {\rm Im}(P_0^s) - {\rm Im}(P_0^{s+1})$ et que $Q_t^sx = 0$, $\forall$ $t=q, \cdots, r$. Alors $r-q \leq d-2$.\qed
\end{Lem}

On laisse la d\'emonstration de ce lemme au lecteur. Voici quelques indications. D'abord, comme dans la section \ref{Se3}, on d\'efinit un $g-$\'echange ($g \geq 0$) entre deux mon\^omes basiques $\alpha$ et $\beta$ de la mani\`ere suivante: un tel $g-$\'echange existe entre $\alpha$ et $\beta$ s'il existe $i$ et $j$ tels que ces deux mon\^omes constituent, \`a un ordre (entre $i$ et $j$) pr\`es, une paire de la forme
$$\alpha = (u_1, \cdots, u_{i-1}, {2u_i+1}, u_{i+1}, \cdots, u_{j-1}, {2u_j+2p^g}, u_{j+1}, \cdots, u_d)$$
$$\beta = (u_1, \cdots, u_{i-1}, {2u_i+2p^g}, u_{i+1}, \cdots, u_{j-1}, {2u_j+1}, u_{j+1}, \cdots, u_d).\leqno{\hbox{et}}$$
On dit aussi que c'est un $g-$\'echange en $i-$\`eme position pour $\alpha$ avec $\beta$. En utilisant cette nouvelle d\'efinition de $g-$\'echange, et les autres d\'efinitions restant inchang\'ees, on aura la m\^eme proposition que la proposition \ref{Pr2} pour le cas $p$ premier impair. Ensuite on ach\`eve la d\'emonstration en construisant une $(q, r)-$classe comme dans la d\'emonstration du corollaire \ref{Co1}.

\begin{Lem}\label{Le6}
Soit $x$ un \'el\'ement de $H^\ast (B(\Bbb{Z}/p)^{\oplus d};$ $\Bbb{Z}/p)$. Si ${\mathcal A}_px$ contient une lacune $(|x|, |x|+l]$, $l \geq 2(p-1)p^k$ pour un certain $k \geq d-3$, alors $x \in {\rm Im}(P_0^{k-d+3})$.
\end{Lem}

\begin{proof}[D\'emonstration]
Soit $\alpha$ tel que $x = P_0^\alpha x' \in {\rm Im}(P_0^\alpha) - {\rm Im}(P_0^{\alpha+1})$. Si $\alpha > k+1$, on a $\alpha \geq k-d+3$. Si $\alpha \leq k+1$, alors l'existence de la lacune dans ${\mathcal A}_px$ implique que $\forall$ $t=0, \cdots, k-\alpha$,
$$P_0^\alpha P^{p^t}x' = P^{p^{t+\alpha}}P_0^\alpha x' = P^{p^{t+\alpha}}x = 0,$$
et que
$$P_0^\alpha \beta x' = P^{p^{\alpha-1}}P_0^\alpha x' = P^{p^{\alpha-1}}x = 0.$$
Donc $\beta x' = 0$ et $P^{p^t}x' = 0$, $\forall$ $t=0, \cdots, k-\alpha$. D'o\`u, $Q_{t'}x' = 0$, $\forall$ $t'=0, \cdots, k-\alpha+1$. D'apr\`es le lemme \ref{Le5}, on a $k-\alpha+1 \leq d-2$, d'o\`u $\alpha \geq k-d+3$.
\end{proof}

\begin{Thm}[\cite{Li62}, \cite{SY61}]\label{Th5} Soit $X$ un espace ou un spectre, $k \geq 1$, soit $x \in H^n(X; \Bbb{Z}/p)$ tel que $\beta x = 0$ et $P^{p^i} x = 0$, $\forall$ $i<k$, alors $P^{p^k} x \in \sum_{i<k} {\rm Im}(P^{p^i}) + {\rm Im}(\beta)$.
\end{Thm}

\begin{Def}
Soit $M$ un module instable infini connexe dont la filtration nilpotente est de longueur finie. Les quotients $nil_sM / nil_{s+1}M$ non-triviaux s'\'ecrivent sous la forme $\Sigma^{m_i} R_{m_i}$, $R_{m_i}$ r\'eduits, $i=1, \cdots, t$, $m_1 < \cdots < m_t$. Notons que l'un au moins des $R_{m_i}$ est infini. Supposons qu'il existe des entiers $d$ et $\alpha_d$ tels que tous les $R_{m_i}$ se plongent dans $H^\ast (B(\Bbb{Z}/p)^{\oplus d}; \Bbb{Z}/p)^{\oplus \alpha_d}$. Notons $I \subset \{ 1, \cdots, t \}$ le sous-ensemble des $i$ tels que $R_{m_i}$ soit infini, et soit $n_{1,i} < n_{2,i} < \cdots$ les degr\'es en lesquels ce module est non-trivial.

Soit $\delta$ tel que $p^\delta \geq t > p^{\delta-1}$. Le module instable $M$ sera dit {\sl de type $\mathcal T$} s'il contient une lacune $(s, s+l]$ avec $s \geq \min\{ m_i+n_{1,i}\ |\ i \in I \}$ et
$$\begin{array}{rcl}
l& \geq& \max\{ 2(m_t+1)(p-1)p^{d+2}, n_{j+1,i}-n_{j,i}\ |\ i \in I, \\
& & \qquad \qquad \qquad j=1, \cdots, 1+(d+\delta)(p-1)^2p^{d-2} \}.
\end{array}$$
\end{Def}

\begin{Cond}\label{Cond2}
Soit $M$ un module instable connexe dont la filtration nilpotente est de longueur finie. En utilisant les notations introduites dans la d\'efinition pr\'ec\'edente, on dira que $M$ v\'erifie la condition \ref{Cond2} si
$$m_j-m_i \neq 1,2(p-1), \qquad \forall\ 1 \leq i, j \leq t.$$
\end{Cond}

\begin{Thm}\label{Th6}
Soit $M$ un module qui est une suspension it\'er\'ee (positive ou n\'egative) d'un module instable connexe dont la filtration nilpotente est de longueur finie, qui est de type $\mathcal T$ et v\'erifie la condition \ref{Cond2}. Alors $M$ n'est pas r\'ealisable, i.e., il n'existe aucun espace ou spectre $X$ tel que $M = \tilde{H}^\ast(X; \Bbb{Z}/p)$.
\end{Thm}

{\bf D\'emonstration}\qua L'id\'ee essentielle de la d\'emonstration de ce th\'eor\`eme est la m\^eme que celle de la d\'emonstration du th\'eor\`eme \ref{Th4}. N\'eanmoins, certains aspects du cas d'un nombre premier impair apparaissent, non seulement on utilise le th\'eor\`eme \ref{Th5} au lieu du th\'eor\`eme d'Adams, mais aussi on a besoin de reconstituer les calculs pour le cas d'un nombre premier impair. On donne dans la suite une esquisse de la d\'emonstration de ce th\'eor\`eme, afin d'illustrer certains changements n\'ecessaires par rapport \`a celle du th\'eor\`eme \ref{Th4}.

On note, comme dans la d\'emonstration du th\'eor\`eme \ref{Th4}, qu'il suffit de prouver l'\'enonc\'e pour les modules instables connexes dont la filtration nilpotente est de longueur finie, qui est de type $\mathcal T$ et v\'erifie la condition \ref{Cond2}.

Ensuite, on raisonne par l'absurde. Supposons qu'il existe un tel module instable connexe $M$ qui est la cohomologie r\'eduite d'un espace ou d'un spectre, et que la lacune $(n, n+l]$ soit la premi\`ere dans $M$, avec $n \geq \min\{ m_i+n_{1,i}\ |\ i \in I \}$, de longueur
$$\begin{array}{rcl}
l& \geq& \text{max}\{ 2(m_t+1)(p-1)p^{d+2}, n_{j+1,i}-n_{j,i}\ |\ i \in I,\\
& & \qquad \qquad \qquad j=1, \cdots, 1+(d+\delta)(p-1)^2p^{d-2} \},
\end{array}$$
o\`u $\delta$ est l'unique entier tel que $p^\delta \geq t > p^{\delta-1}$. Soit $k$ l'unique entier ($\geq d+2$) tel que $2(p-1)p^{k+1} > l \geq 2(p-1)p^k$.

Si on note $\forall$ $i \in I$,
$$(n-m_i \geq)\ r_{0,i} > r_{1,i} > \cdots$$
les degr\'es inf\'erieurs \`a $n+l-m_i$, pour lesquels $R_{m_i}$ est non-trivial. Alors si on peut montrer que pour tout \'el\'ement $x$ de $M$ de degr\'e inf\'erieur ou \'egal \`a $n$, ${\mathcal A}_px$ contient une lacune $(|x|, n+l]$, on peut montrer par r\'ecurrence que
$$\forall\ w \geq j(p-1)^2p^{d-2},\qquad 2p^{k-d+j+2} | r_{w,i}.$$
Puisque pour tout \'el\'ement $x$ de $M$, et donc pour sa r\'eduction $\sigma^{m_i} x'$ dans $\Sigma^{m_i} R_{m_i}$, de degr\'e inf\'erieur ou \'egal \`a $n$, ${\mathcal A}_px$ contient une lacune $(|x|, n+l]$ et ${\mathcal A}_px'$ contient une lacune $(|x'|, n+l-m_i]$, le lemme \ref{Le6} montre que le cas $j = 0$ est vrai.

Supposons que c'est vrai pour $j$, alors
$$\begin{array}{rl}
&      (n+l-m_i)-r_{(j+1)(p-1)^2p^{d-2},i}\\
\geq&  l+r_{0,i}-r_{(j+1)(p-1)^2p^{d-2},i}\\
=&     l+\sum_{w=0}^{(j+1)(p-1)^2p^{d-2}-1} (r_{w,i}-r_{w+1,i})\\
=&     l+\sum_{h=0}^{j} \sum_{w=h(p-1)^2p^{d-2}}^{(h+1)(p-1)^2p^{d-2}-1} (r_{w,i}-r_{w+1,i})\\
\geq&  2(p-1)p^k+\sum_{h=0}^{j} (p-1)^2p^{d-2} \cdot 2p^{k-d+h+2}\\
=&     2(p-1)p^k+2(p-1)^2\sum_{h=0}^{j} p^{k+h}\\
=&     2(p-1)p^k+2(p-1)(p^{k+j+1}-p^k)\\
=&     2(p-1)p^{k+j+1}.
\end{array}$$
Puisque pour les \'el\'ements $x'$ en degr\'e $r_{w,i}$, ${\mathcal A}_px'$ contient une lacune $(r_{w,i}, n+l-m_i]$, le lemme \ref{Le6} montre que $2p^{k-d+j+3} | r_{w,i}$, $\forall\ w \geq (j+1)(p-1)^2p^{d-2}$.

Rappelons que par hypoth\`ese, il y a bien (au moins) $(d+\delta)(p-1)^2p^{d-2} + 1$ valeurs pour l'indice $i$ de $r_i$. On peut donc poser $w_0 = (d+\delta)(p-1)^2p^{d-2}$, alors $2p^{k+\delta+2} | r_{w,i}$,$\forall$ $w \geq w_0$. Donc en degr\'e plus petit que $r_{w_0,i}$, $R_{m_i}$ ($i \in I$) ne contient que des lacunes de longueur plus grande que ou \'egale \`a $2p^{k+\delta+2}-1$. Donc on peut choisir un $d_0 \leq n$ tel qu'en degr\'e inf\'erieur ou \'egal \`a $d_0$, les $\Sigma^{m_i} R_{m_i}$ ($i \in I$) ne contiennent que des lacunes de longueur plus grande que ou \'egale \`a $2p^{k+\delta+2}-1$, et on suppose de plus qu'il en existe au moins une en degr\'e inf\'erieur ou \'egal \`a $d_0$.

On peut donc choisir, comme \`a la fin de la d\'emonstration du th\'eor\`eme \ref{Th4}, une lacune en degr\'e inf\'erieur ou \'egal \`a $d_0$, dans l'un des $\Sigma^{m_i} R_{m_i}$ ($i \in I$), telle qu'elle contient au plus $t-1$ degr\'es en lesquels $M$ est non-trivial. Donc il existe une lacune $(n', n'+l']$, avec $n' \geq \min\{ m_i+n_{1,i}\ |\ i \in I \}$ et $n'+l' \leq d_0$ ($\leq n$), de longueur
$$\begin{array}{rcl}
l'& \geq& \frac{1}{t} \cdot 2p^{k+\delta+2} - 1 \\
& \geq& 2p^{k+2} - 1 \\
& >& 2(p-1)p^{k+1} \\
& >& l \\
& \geq& \max\{ 2(m_t+1)(p-1)p^{d+2}, n_{j+1,i}-n_{j,i}\ |\ i \in I, \\
& & \qquad \qquad \qquad j=1, \cdots, 1+(d+\delta)(p-1)^2p^{d-2} \}.
\end{array}$$
Ceci est contradictoire au choix de $(n, n+l]$.

Pour finir la d\'emonstration du th\'eor\`eme, il reste donc \`a montrer le

\begin{Lem}\label{Prop04}
Pour tout $x$ tel que $|x| \leq n$, le module ${\mathcal A}_px$ contient la lacune $(|x|, n+l]$.
\end{Lem}

{\bf D\'emonstration}\qua Comme dans la d\'emonstration du lemme \ref{Prop03}, on raisonne par l'absurde. A tout \'el\'ement $x \in M$, on associe son degr\'e de nilpotence, c'est-\`a-dire, l'entier $m_x$ tel que $x \in nil_{m_x}M - nil_{m_x+1}M$.

Soit $x$ un \'el\'ement non nul de degr\'e maximal tel que ${\mathcal A}_px$ n'est pas r\'eduit \`a $\{0\}$ dans l'intervalle $(|x|, n]$. Soit donc $y \in {\mathcal A}_px$ de degr\'e minimal tel que $y \neq 0$, $|x| < |y| \leq n$, $\bar{y} \in \Sigma^{m_y} R_{m_y}$ sa r\'eduction que l'on note $\sigma^{m_y} v$, $v \in R_{m_y}$. On a, \`a l'aide du lemme \ref{Le6}, $v \in {\rm Im}(P_0^{k-d+3})$. On sait donc que $|y|-m_y$ est divisible par $2p^{k-d+2}$ et on \'ecrit $|y|-m_y =: 2p^{k-d+2}l_y$.

Soit de m\^eme la r\'eduction $\bar{x} \in \Sigma^{m_x} R_{m_x}$. Notons $\bar{x} = \sigma^{m_x} u$, $u \in {\rm Im}(P_0^s) - {\rm Im}(P_0^{s+1})$.

\begin{Lem}
Pour $s \leq k-d+1$, $t \geq 1$, $2p^{k-d+2}$ ne divise pas le degr\'e de $Q_t^s u$.
\end{Lem}

\begin{proof}[D\'emonstration]
En effet on a
$$|Q_t^s u| = \left\{ \begin{array}{ll}
    |u|+2(p^{t+s}-p^{s-1})& s \geq 1\\
    |u|+(2p^t-1)& s=0
    \end{array} \right.$$
Consid\'erons le module engendr\'e par $u$ dans $R_{m_x}$. Si ${\mathcal A}_p u$ contient la lacune $(|u|, n+l-m_x]$, le degr\'e de $u$ est divisible par $2p^{k-d+2}$, et dans ce cas $2p^{k-d+2}$ ne divise pas $|Q_t^s u|$ pour $s \leq k-d+1$.

Supposons que ${\mathcal A}_p u$ ne contienne pas la lacune $(|u|, n+l-m_x]$ et soit $t = Pu$ l'\'el\'ement non nul du plus bas degr\'e avec $|t| \leq n-m_x$. On sait que le degr\'e de l'op\'eration $P$ est $2(p-1)p^\gamma$ ou 1. Pour la m\^eme raison que dans la d\'emonstration du lemme \ref{Le3}, on a $\gamma \geq s$. Comme ${\mathcal A}_p t$ contient la lacune $(|t|, n+l-m_x]$, $|t|$ est divisible par $2p^{k-d+2}$ et peut donc \^etre \'ecrit de la forme $|t| = 2p^{k-d+2}q$. Alors
$$|Q_t^s u| = \left\{ \begin{array}{ll}
    2p^{k-d+2}q - 2(p-1)p^\gamma + 2(p^{t+s} - p^{s-1})& s \geq 1\\
    2p^{k-d+2}q - 2(p-1)p^\gamma + (2p^t-1) \ \text{ou}& \\
    2p^{k-d+2}q - 1 + (2p^t-1)& s = 0
    \end{array} \right.$$
et $2p^{k-d+2}$ ne divise pas $|Q_t^s u|$.
\end{proof}

Supposons d'abord que $s \leq k-d+1$ et que $t \geq 1$. Tant que $|Q_t^s x| \leq n+l$, on a n\'ecessairement $Q_t^s u = 0$ pour des raisons de degr\'e. Or pour $1 \leq t \leq k-s$,
$$\begin{array}{rcl}
|Q_t^s x|& =& \left\{ \begin{array}{ll}
    |x|+2(p^{t+s}-p^{s-1})& s \geq 1\\
    |x|+(2p^t-1)& s = 0
    \end{array} \right.\\
&\leq& \left\{ \begin{array}{ll}
    n+2(p^k-p^{s-1})& s \geq 1\\
    n+(2p^k-1)& s = 0
    \end{array} \right.\\
&<& n+2p^k\\
&<& n+2(p-1)p^k\\
&\leq& n+l.
\end{array}$$
Donc, d'apr\`es le lemme \ref{Le5},
$$k-s-1 \leq d-2, \quad \text{soit}\ s \geq k-d+1$$
et donc $|x|-m_x =: 2p^{k-d}l_x$.

Maintenant on peut trouver une contradiction comme dans la d\'emonstration du lemme \ref{Prop03}. A l'aide du th\'eor\`eme \ref{Th5} et de la base multiplicative de ${\mathcal A}_p$, on sait qu'en degr\'e sup\'erieur ou \'egal \`a $|x|+1$, l'\'el\'ement non nul du plus bas degr\'e dans ${\mathcal A}_px$ ne peut \^etre que $y = \beta x$ ou $y = P^1 x$. On a alors $|x|+1 = |y|$ ou $|x|+2(p-1) = |y|$. Donc
$$\begin{array}{rrcl}
& m_x+2p^{k-d}l_x+1& =& m_y+2p^{k-d+2}l_y,\\
\text{ou}& m_x+2p^{k-d}l_x+2(p-1)& =& m_y+2p^{k-d+2}l_y.
\end{array}$$
Donc, $m_y-m_x = 1$ ou $2(p-1)$ mod $2p^{k-d}$.

Puisque $2(p-1)p^{k+1} > l \geq 2(m_t+1)(p-1)p^{d+2}$, donc
$$\begin{array}{rcl}
2p^{k-d}& >& 2p(m_t+1)\\
& \geq& 6m_t+2p\\
& >& \left\{ \begin{array}{l}
    m_y+m_x+2(p-1)\\
    m_y+m_x+1
    \end{array} \right.\\
& \geq& \left\{ \begin{array}{l}
    |m_y-m_x-2(p-1)|\\
    |m_y-m_x-1|
    \end{array} \right.
\end{array}$$
Donc $m_y-m_x = 1$ ou $2(p-1)$. Or cette \'egalit\'e n'a pas de solution \`a cause de la condition \ref{Cond2}, l'existence d'un tel $x$ est donc contradictoire. Donc pour tout \'el\'ement $x$ de $M$ de degr\'e inf\'erieur \`a $n$, ${\mathcal A}_px$ contient une lacune $(|x|, n+l]$.

{\bf Fin de la d\'emonstration du Lemme \ref{Prop04}} \qed

{\bf Fin de la d\'emonstration du Th\'eor\`eme \ref{Th6}}\qed

\section*{Appendice: Notes sur les op\'erations $Q_t^s$}
\addcontentsline{toc}{section}{Appendice: Notes sur les op\'erations $Q_t^s$}
\small

Dans la base de Milnor de l'alg\`ebre de Steenrod ${\mathcal A}_2$, on a des op\'erations $Q_t$, $t \geq 0$, d\'efinies r\'ecursivement par les relations suivantes:

(1)\qua $Q_0 = Sq^1$;

(2)\qua $Q_{t+1} = [Sq^{2^{t+1}}, Q_t]$.

Ces op\'erations ont des bonnes propri\'et\'es, plus pr\'ecis\'ement (Milnor, \cite{Mi58}, $\S 6$),

\begin{Prop}\label{Pr6}
Soient $M, N$ deux modules instables, et soient $l, r, t \geq 0$.

(1)\qua $\forall$ $x \in M$,
$Q_rQ_t x = Q_tQ_r x$ et $Q_t^2 x = 0$.

(2)\qua $\forall$ $x \in M$ et $y \in N$,
$Q_t(x \otimes y) = Q_tx \otimes y + x \otimes Q_ty$.

(3)\qua Soit $u$ le g\'en\'erateur de $H^\ast (B(\Bbb{Z}/2); \Bbb{Z}/2)$ en degr\'e 1,

$Q_tu^{2l} = 0$ et $Q_tu^{2l+1} = u^{2l+2^{t+1}}$.
\end{Prop}

Inspir\'e par ces propri\'et\'es, on a construit dans la section 2 les op\'erations $Q_t^s$, $s, t \geq 0$, qui poss\`edent aussi ces propri\'et\'es sur ${\rm Im}(Sq_0^s)$. Dans le reste de cet appendice, on d\'ecrit quelques propri\'et\'es \'el\'ementaires de ces op\'erations. Puis, on donne \`a la fin les d\'emonstrations du lemme \ref{Le1} et du corollaire \ref{Pr1}.

D'abord, on compare ces op\'erations avec les op\'erations connues, $P_{t+1}^s$, $s, t \geq 0$. Voici deux propri\'et\'es \'el\'ementaires:

(1)\qua $Q_t^0 = Q_t = P_{t+1}^0$, $Q_0^s = Sq^{2^s} = P_1^s$.

(2)\qua $Q_t^s$ est une op\'eration de degr\'e $2^s(2^{t+1}-1)$.

D'apr\`es ces deux propri\'et\'es, plus le fait que $P_{t+1}^s$ est aussi une op\'eration de degr\'e $2^s(2^{t+1}-1)$, une question curieuse est de savoir quand les deux op\'erations $Q_t^s$ et $P_{t+1}^s$ coincident. En fait, quand $st \neq 0$, il semble que le seul cas o\`u ces deux op\'erations coincident est $Q_1^1 = P_2^1$.

Pour effectuer le calcul des $Q_t^s$, on note que l'on a une autre fa\c{c}on de d\'efinir les op\'erations $Q_t^s$, $s, t \geq 0$:

(1)\qua $Q_0^s = P_1^s$;

(2)\qua $Q_{t+1}^s = [P_1^{t+s+1}, Q_t^s]$.

Avec cette d\'efinition, on peut exprimer $Q_t^s$ en terme de la base de Milnor, \`a l'aide de la formule multiplicative de cette base. Voici quelques calculs qui comparent les $Q_t^s$ et $P_{t+1}^s$:

(1)\qua $Q_1^1 = P_2^1$, $Q_1^2=P_2^2+Sq(3,3)$, $Q_1^3=P_2^3+Sq(6,6)+Sq(3,7)$.

(2)\qua $Q_2^1=P_3^1+Sq(7,0,1)+Sq(4,1,1)$.

Pour finir cet appendice, on donne ici les d\'emonstrations du lemme \ref{Le1} et du corollaire \ref{Pr1}.

{\bf D\'emonstration du Lemme \ref{Le1}}\qua
Soient $M$ un module instable et $x$ un \'el\'ement de $M$.

Quand $t=0$,
$$Q_0^{s+r}Sq_0^s x = Sq^{2^{s+r}}Sq_0^s x = Sq_0^sSq^{2^r} x = Sq_0^sQ_0^r x.$$

Si on suppose que $Q_t^{s+r}Sq_0^s x = Sq_0^sQ_t^r x$ pour $t \leq t_0$, alors quand $t = t_0+1$, on a
$$\begin{array}{lll}
Q_{t_0+1}^{s+r}Sq_0^s x&=&   (Sq^{2^{t_0+s+r+1}}Q_{t_0}^{s+r} - Q_{t_0}^{s+r}Sq^{2^{t_0+s+r+1}})Sq_0^s x\\
&=&     Sq^{2^{t_0+s+r+1}}Q_{t_0}^{s+r}Sq_0^s x - Q_{t_0}^{s+r}Sq^{2^{t_0+s+r+1}}Sq_0^s x\\
&=&     Sq^{2^{t_0+s+r+1}}Sq_0^sQ_{t_0}^r x - Q_{t_0}^{s+r}Sq_0^sSq^{2^{t_0+r+1}} x\\
&=&     Sq_0^sSq^{2^{t_0+r+1}}Q_{t_0}^r x - Sq_0^sQ_{t_0}^rSq^{2^{t_0+r+1}} x\\
&=&     Sq_0^s(Sq^{2^{t_0+r+1}}Q_{t_0}^r - Q_{t_0}^rSq^{2^{t_0+r+1}}) x\\
&=&     Sq_0^sQ_{t_0+1}^r x.
\end{array}$$
Donc, par r\'ecurrence (sur $t$), on a $Q_t^{s+r}Sq_0^s x = Sq_0^sQ_t^r x$, $\forall$ $r, s, t$.\endproof

{\bf D\'emonstration du Corollaire \ref{Pr1}}\qua
Les propri\'et\'es de $Q_t$ utilis\'ees ci-dessous sont dans la proposition \ref{Pr6}.

(1)\qua Soit $M$ un module instable. Par le lemme \ref{Le1}, on a $\forall$ $x \in M$,
$$\begin{array}{rllllll}
Q_r^sQ_t^sSq_0^s x&=&        Q_r^sSq_0^sQ_t x&=&        Sq_0^sQ_rQ_t x&&\\
&=&     Sq_0^sQ_tQ_r x&=&        Q_t^sSq_0^sQ_r x&&\\
&=&     Q_t^sQ_r^sSq_0^s x,&&&&\\
\text{ et }\quad (Q_t^s)^2Sq_0^s x&=&      Q_t^sSq_0^sQ_t x&=&        Sq_0^s(Q_t)^2 x&=& 0.
\end{array}$$

(2)\qua Supposons qu'il existe $x' \in M$, $y' \in N$ tels que $x = Sq_0^s(x')$, $y = Sq_0^s(y')$, alors
$$\begin{array}{lll}
Q_t^s(x \otimes y)&=&   Q_t^s(Sq_0^s(x') \otimes Sq_0^s(y'))\\
&=&     Q_t^sSq_0^s(x' \otimes y')\\
&=&     Sq_0^sQ_t(x' \otimes y')\\
&=&     Sq_0^s(Q_tx' \otimes y' + x' \otimes Q_ty')\\
&=&     Sq_0^sQ_tx' \otimes Sq_0^sy' + Sq_0^sx' \otimes Sq_0^sQ_ty'\\
&=&     Q_t^sx \otimes y + x \otimes Q_t^sy.
\end{array}$$

(3)\qua On a 
$$\begin{array}{rllllll}
Q_t^sSq_0^su^{2l}&=&       Sq_0^sQ_tu^{2l}&=& 0,&&\\
Q_t^sSq_0^su^{2l+1}&=&     Sq_0^sQ_tu^{2l+1}&=&       Sq_0^su^{2l+2^{t+1}}&=&    u^{2^s(2l+2^{t+1})}.
\end{array}\eqno{\qed}$$\normalsize

\renewcommand{\refname}{Bibliographie}

\Addresses\recd

\end{document}

%% file: agtout.tex
\def\ifplaintex{\expandafter\ifx\csname documentclass\endcsname\relax}

\def\gtp{{\mathsurround=0pt\it $\cal G\mskip-2mu$eometry \&\ 
$\cal T\!\!$opology $\cal P\!$ublications}}  % GT publications

\def\recd{{\small Received:\qua\receiveddate\ifx\reviseddate\relax
\else\qquad Revised:\qua\reviseddate\fi\par}} 

\def\lognumber#1{\def\thelognumber{#1}}
\def\volumenumber#1{\def\thevolumenumber{#1}}
\def\volumeyear#1{\def\thevolumeyear{#1}}
\def\papernumber#1{\def\thepapernumber{#1}}
\def\pagenumbers#1#2{\def\startpage{#1}\def\finishpage{#2}}
\def\published#1{\def\publishdate{#1}}

\def\received#1{\def\receiveddate{#1}}
\def\revised#1{\def\reviseddate{#1}}
\def\accepted#1{\def\accepteddate{#1}}
\def\asciititle#1{\def\theasciititle{#1}}
\def\covertitle#1{\def\thecovertitle{#1}}

\def\asciiaddress#1{\def\theasciiaddress{#1}}

\long\def\asciiabstract#1{\long\def\theasciiabstract{#1}}
\def\asciikeywords#1{\def\theasciikeywords{#1}}

\let\\\par\let\thelognumber\relax\let\thevolumenumber\relax
\let\thepapernumber\relax\let\thevolumeyear\relax\let\startpage\relax
\let\finishpage\relax\let\publishdate\relax\let\receiveddate\relax
\let\reviseddate\relax\let\accepteddate\relax\let\theasciititle\relax
\let\thecovertitle\relax\let\theasciiauthors\relax\let\theasciiaddress\relax
\let\theasciiabstract\relax\let\theasciikeywords\relax

\let\theasciiemail\relax

\ifplaintex
\font\logobig=cmssbx10 scaled 3836
\font\logomed=cmssbx10 scaled 2557
\else
\font\logobig=cmssbx10 scaled 4200
\font\logomed=cmssbx10 scaled 2800
\fi

\long\def\makeagttitle{   %%% start of definition of \makeagttitle
\count0=\startpage
\agt\hfill      %   Journal title (top left) 
\hbox to 45truept{\vbox to 0pt{\vglue -13truept{\logomed A\kern -.37em{\logobig 
T}\kern -.38em G}\vss}\hss}
\break
{\small Volume \thevolumenumber\ (\thevolumeyear)
\startpage--\finishpage\nl
Published: \publishdate}

\vglue .25truein

{\parskip=0pt\leftskip 0pt plus
1fil\def\\{\par\smallskip}{\Large\bf\thetitle}\par\medskip} \vglue
0.05truein

{\parskip=0pt\leftskip 0pt plus 1fil\def\\{\par}{\sc\theauthors}
\par\medskip}%
 
\vglue 0.03truein 

{\small\leftskip 25truept\rightskip 25truept{\bf Abstract}\stdspace\theabstract

{\bf AMS Classification}\stdspace\theprimaryclass
\ifx\thesecondaryclass\relax\else; \thesecondaryclass\fi\par
{\bf Keywords}\stdspace \thekeywords\par}\vglue 7truept

}   %%%% end of definition of \makeagttitle

\ifplaintex
\font\phead=cmsl9 scaled 950
\font\pnum=cmbx10 scaled 913
\font\pfoot=cmsl9 scaled 950
\headline{\vbox to 0pt{\vskip -4.5mm\line{\small\phead\ifnum
\count0=\startpage ISSN 1472-2739 (on-line) 1472-2747 (printed)
\hfill {\pnum\folio}\else\ifodd\count0\def\\{ }% 
\ifx\theshorttitle\relax\thetitle\else\theshorttitle\fi\hfill{\pnum\folio}
\else\def\\{ and }{\pnum\folio}\hfill\ifx\theshortauthors\relax\theauthors
\else\theshortauthors\fi\fi\fi}\vss}}
\footline{\vbox to 0pt{\vglue 0mm\line{\small\pfoot\ifnum\count0=\startpage
\copyright\ \gtp\hfill\else
\agt, Volume \thevolumenumber\ (\thevolumeyear)\hfill\fi}\vss}}
\else
\font=cmsl9 scaled 1050
\font\lnum=cmbx10 
\font=cmsl9 scaled 1050
\makeatletter
\def\@oddhead{{\small\lhead\ifnum\count0=\startpage ISSN 1472-2739 
(on-line) 1472-2747 (printed)\hfill {\lnum\number\count0}\else\ifodd\count0
\def\\{ }\ifx\theshorttitle\relax \thetitle \else\theshorttitle\fi\hfill
{\lnum\number\count0}\else\def\\{ and }{\lnum\number\count0}
\hfill\ifx\theshortauthors\relax 
\theauthors\else\theshortauthors\fi\fi\fi}}\def\@evenhead{\@oddhead}
\def\@oddfoot{\small\lfoot\ifnum\count0=\startpage\copyright\ \gtp\hfill\else
\agt, Volume \thevolumenumber\ (\thevolumeyear)\hfill\fi}
\def\@evenfoot{\@oddfoot}
\makeatother
\fi
\let\maketitlepage\makeagttitle
\let\makeshorttitle\maketitlepage
\let\maketitle\maketitlepage

\newwrite\gtoutfile
\long\gdef\makeheadfile{  %%% start of definition of \makeheadfile
{\def\\{, }\def\s{ }
\immediate\openout\gtoutfile head.xxx
\immediate\write\gtoutfile{Proxy-for: \ifx\theasciiauthors\relax
\theauthors\else\theasciiauthors\fi\s<\ifx\theasciiemail\relax\theemail\else\theasciiemail\fi>}
\immediate\write\gtoutfile{\noexpand\\}
\immediate\write\gtoutfile{Authors: \ifx\theasciiauthors\relax
\theauthors\else\theasciiauthors\fi}
{\def\\{ }\immediate\write\gtoutfile{Title: \ifx\theasciititle\relax
\thetitle\else\theasciititle\fi}}
\immediate\write\gtoutfile{Subj-class: GT or SG, GR etc}
\immediate\write\gtoutfile{MSC-class: \theprimaryclass\ifx\thesecondaryclass\relax\else, \thesecondaryclass\fi}
\immediate\write\gtoutfile{Journal-ref: Algebraic and Geometric Topology \thevolumenumber\s
(\thevolumeyear) \startpage-\finishpage}
\immediate\write\gtoutfile{Comments: Published by Algebraic and
Geometric Topology at}
\immediate\write\gtoutfile{\s\s\s  http://www.maths.warwick.ac.uk/agt/AGTVol\thevolumenumber/agt-\thevolumenumber-\thepapernumber.abs.html}
\immediate\write\gtoutfile{\noexpand\\}
\immediate\write\gtoutfile{}
\ifx\theasciiabstract\relax
\immediate\write\gtoutfile{\theabstract}\else
\immediate\write\gtoutfile{\theasciiabstract}\fi
\immediate\write\gtoutfile{}
\immediate\write\gtoutfile{\noexpand\\}
\immediate\write\gtoutfile{}
\immediate\closeout\gtoutfile}}  %%% end of definition of \makeheadfile

\def\maketitlepage{\makeagttitle\makeheadfile}
\let\makeshorttitle\maketitlepage
\let\maketitle\maketitlepage